\documentclass[english]{article}
\usepackage[latin9]{inputenc}
\usepackage{textcomp}
\usepackage{amsmath}
\usepackage{amssymb}
\usepackage{stackrel}

\makeatletter
\newcommand{\lyxaddress}[1]{
\par {\raggedright #1
\vspace{1.4em}
\noindent\par}
}

\makeatother

\usepackage{babel}
\begin{document}

\title{Kernel of Trace Operator of Sobolev Spaces on Lipschitz Domain}

\author{I-Shing Hu}
\date{}
\maketitle
\begin{abstract}
We are going to show that on bounded Lipschitz domain $D$: both $C_{c}^{\infty}(D)$,
the set of smooth functions on $D$ with compact support, and $C_{0}^{\infty}(D)$,
the set of smooth functions on $D$ with (extension) zero boundary,
are dense in $W^{1,p}\left(D\right)$, $p\in[1,\infty)$. A proof
can be found in Ne\v{c}as's monograph \cite{key-2}, Theorem 4.10,
§2.4.3. 

Our main result in this note is that: we find another proof by showing
that both closures is the same as kernel of trace operator 
\[
T:\,W^{1,p}(D)\rightarrow L^{p}(\partial D).
\]

via some change of variables formulas from Evans and Gariepy's textbook
\cite{key-4} for Lipschitz coordinate transformation, to extend the
proof of Theorem 2 in §5.5 of Evans' widespread PDE textbook \cite{key-3},
from $\mathcal{C}^{1}$ to Lipschitz domain.
\end{abstract}
Firstly we review some analysis on Lipschitz domain. Let $C_{c}^{\infty}(D)$
be the set of smooth functions on $D$ with compact support, and $C_{0}^{\infty}(D)$
be the set of smooth functions on $D$ with (extension) zero boundary.
It is worth to note that in Grisvard's classic \cite{key-1}, Corollary
1.5.1.6 in §1.5.1, states without proof a much more general result
which covers above.

\section{Preliminary: Lipschitz Domain}

\subsection{Lipschitz homeomorphism }

Let $D\subset\mathbb{R}^{d}$ be a Lipschitz domain. First we state
the usual coordinate transformation from Appendix C.1 in \cite{key-3}.
Fix $a\in\partial D$, there exist $r>0$ and Lipschitz map $\gamma:\,\mathbb{R}^{d-1}\rightarrow\mathbb{R}$
such that (upon relabeling and reorienting coordinates axes if necessary)
\[
D\cap B(a,r)=\left\{ x\in B(a,r)\mid x_{d}>\gamma(x_{1},\ldots,\,x_{d-1})\right\} .
\]
Let $F:\,\mathbb{R}^{d}\rightarrow\mathbb{R}^{d}$ be Lipschitz coordinate
transformation 
\[
F(x_{1},\ldots,x_{d-1},x_{d}):=(x_{1},\ldots,x_{d-1},\,x_{d}-\gamma(x_{1},\ldots,x_{d-1})),
\]
and obviously it has inverse $F^{-1}$ as
\[
F^{-1}(y_{1},\ldots,y_{d-1},y_{d})=(y_{1},\ldots,y_{d-1},\,y_{d}+\gamma(y_{1},\ldots,y_{d-1})),
\]
which is also continuous. Hence $D$ is homeomorphic to an open subset
of $\mathbb{R}_{+}^{d}:=\left\{ (y_{1},\ldots,y_{d-1},y_{d})\in\mathbb{R}^{d}\mid y_{d}>0\right\} $,
and obviously $\partial D$ is homeomorphic to a closed subset of
$\left\{ (y_{1},\ldots,y_{d-1},y_{d})\in\mathbb{R}^{d}\mid y_{d}=0\right\} $.
To calculate functions on (an open subset of) $\mathbb{R}_{+}^{d}$
instead of $D$, we let $\left\{ f_{i}:\,U_{i}\subset D\rightarrow\mathbb{R}_{+}^{d}\right\} $
be an atlas of Lipschitz coordinate maps on $D$ as above with $\overline{U}_{i}$
compact, and write $\left\{ \rho_{i}:\,D\rightarrow[0,1]\right\} _{i}$
for a partition of unity associated to $\left\{ U_{i}\right\} _{i}$.
Then $\forall\phi\in W^{1,p}\left(D\right)$, we can map $\phi$ as
\[
\phi\mapsto\underset{i}{\sum}\rho_{i}\phi\circ f_{i}\in W^{1,p}\left(\mathbb{R}_{+}^{d}\right).
\]
It is not difficult to show that this form is invariant under partitions
of unity. For any $a\in\partial D$ and $r>0$, let $\left\{ B(a,r)\mid a\in\partial D\right\} _{r>0}$
be an open cover, then there exists finite sub-cover $\left\{ B_{i}\right\} _{i=1}^{\nu}$.
Then let $D_{i}:=D\cap B_{i}$, $D_{0}:=D-\left(\stackrel[i=1]{\nu}{\cup}D_{i}\right)$,
$\varphi_{i}\in C_{c}^{\infty}(D_{i})$ be a partition of unity on
$D$, i.e.,
\[
1_{D}=\stackrel[i=0]{\nu}{\sum}\varphi_{i}.
\]
Since our estimate is local (on $D_{i}$), we are going to show local
estimate can be extended to global (on $D$). Let $f\in C_{0}^{\infty}(D)$,
$h_{i}\in C_{0}^{\infty}(D_{i})$, then
\[
||f-h||_{W^{1,2}(D)}^{2}\leq c_{1}\stackrel[i=0]{\nu}{\sum}\int\varphi_{i}^{2}\left(|f-h_{i}|^{2}+|Df-Dh_{i}|^{2}\right)+c_{2}\stackrel[i=0]{\nu}{\sum}\int|f-h_{i}|^{2}\left(D\varphi_{i}\right)^{2}.
\]
The only remaining obstacle to prove 
\[
\overline{C_{0}^{\infty}\left(D\right)}=\overline{C_{c}^{\infty}\left(D\right)}
\]
in $W^{1,p}\left(D\right)$, is to integrate under change of variables.
Thereby we need some auxiliary change of variables formulas below.
By Rademacher theorem (see §3.1.2 of \cite{key-4}), Jacobi matrix
$DF$ exists a.e. and Jacobian $JF=1$ a.e. 

\subsection{Change of variables formulas}

We state without proofs two auxiliary change of variables formulas.
They are based on area or co-area formula.

\paragraph{Theorem 1 (Theorem 3.9 in \cite{key-4}).}

Let $f:\,\mathbb{R}^{n}\rightarrow\mathbb{R}^{m}$ be Lipschitz, $n\leq m$.
Then for each $g\in L^{1}(\mathbb{R}^{n})$,

\[
\underset{\mathbb{R}^{n}}{\int}g(x)\,Jf(x)\,dx=\underset{\mathbb{R}^{m}}{\int}\left[\underset{x\in f^{-1}\{w\}}{\sum}g(x)\right]d\mathcal{H}^{n}(w).
\]
\\

\paragraph{Theorem 2 (Theorem 3.11 in \cite{key-4}).}

Let $f:\,\mathbb{R}^{n}\rightarrow\mathbb{R}^{m}$ be Lipschitz, $n\geq m$.
Then for each $g\in L^{1}(\mathbb{R}^{n})$,
\begin{enumerate}
\item For a.e. $y\in\mathbb{R}^{m}$ in Lebesgue measure, $g\in L^{1}(f^{-1}\{y\})$
in (n-m) dimensional Hausdorff measure $\mathcal{H}^{n-m}$, and
\item 
\[
\underset{\mathbb{R}^{n}}{\int}g\,Jf\,dx=\underset{\mathbb{R}^{m}}{\int}\left[\underset{f^{-1}\left\{ w\right\} }{\int}g\,d\mathcal{H}^{n-m}\right]dw.
\]
\end{enumerate}

\section{Trace Operator}

Remark that a proof of trace Theorem on Lipschitz domain can be found
in §4.3 of \cite{key-4}. Adopting from §5.5 of \cite{key-3}, we
can show the existence of trace operator on Lipschitz domain via Theorem
1. 

\subsection{Trace operator on Lipschitz domain}

\paragraph{Theorem 3 (Trace Theorem on Lipschitz domain).}

Let $D\subset\mathbb{R}^{d}$ be a Lipschitz domain. Then for each
$p\in[1,\infty),$ there exists a bounded operator, 
\[
T:\,W^{1,p}(D)\rightarrow L^{p}(\partial D)
\]
such that
\begin{enumerate}
\item $Tu=u|_{\partial D}$ if $u\in W^{1,p}(D)\cap\mathcal{C}(\overline{D})$;
and
\item $\Vert Tu\Vert_{L^{p}(\partial D)}\leq C\Vert u\Vert_{W^{1,p}(D)},$
where $C$ is independent of $u$.\\
\end{enumerate}
Proof: Fix $a\in\partial D$. Then there exist $r>0$ and Lipschitz
map $\gamma$ as in §1.1. Let $f(x)=(x,\gamma(x))$, $x\in\mathbb{R}^{d-1}$
and 
\[
g(x)=\mathtt{1}_{f^{-1}(D\cap B(a,\frac{r}{2}))}|u(f(x))|^{p}.
\]
Remark that we need $p<\infty$. Obviously $f$ is injective and Jacobian
(cf. \cite{key-4}, §3.3.4) $Jf=1+|D\gamma|^{2}\leq c$, where $c$
equals to one plus square of the Lipschitz constant of $\gamma$.
Then
\[
\underset{B(a,\frac{r}{2})\cap\partial D}{\int}|u|^{p}d\mathcal{H}^{d-1}=\underset{\mathbb{R}^{d}}{\int}\left[\underset{x\in f^{-1}\{w\}}{\sum}g(x)\right]d\mathcal{H}^{d-1}(w).
\]
By Theorem 1,
\[
\underset{\mathbb{R}^{d}}{\int}\left[\underset{x\in f^{-1}\{w\}}{\sum}g(x)\right]d\mathcal{H}^{d-1}(w)=\underset{\mathbb{R}^{d-1}}{\int}g(x)\,Jf(x)\,dx\leq c\underset{\mathbb{R}^{d-1}}{\int}g(x)dx.
\]
Let $\zeta\in\mathcal{C}_{c}^{\infty}(B(a,r))$ with $0\leq\zeta\leq1$,
$\zeta=1$ on $B(a,\frac{r}{2})$. Since $\mathbb{R}^{d-1}$ can be
seen as $\{(x_{1},\ldots,x_{d})\in\mathbb{R}^{d}\mid x_{d}=0\}$,
\[
\underset{\mathbb{R}^{d-1}}{\int}g(x)dx\leq-\underset{B(a,r)\cap\{x_{d}\geq0\}}{\int}(\zeta|u|^{p})_{x_{d}}\,dx.
\]
Expand
\[
-\underset{B(a,r)\cap\{x_{d}\geq0\}}{\int}(\zeta|u|^{p})_{x_{d}}\,dx=-\underset{B(a,r)\cap\{x_{d}\geq0\}}{\int}|u|^{p}\zeta_{x_{d}}+p|u|^{p-1}sgn(u)u_{x_{d}}\zeta\,dx.
\]
Since $\partial D$ is compact, for any open cover as above, there
exists a finite sub-cover. For such finite sub-cover, $\sup|D\zeta|$
is uniformly bounded on this finite sub-cover. Hence we have 
\[
\underset{\partial D}{\int}|u|^{p}\,d\mathcal{H}^{d-1}\leq C\underset{D}{\int}|u|^{p}+|Du|^{p}\,dx,
\]
where constant $C$ does not depend on $u$. Write $T:\,W^{1,p}(D)\rightarrow L^{p}(\partial D)$
as
\[
Tu:=u|_{\partial D},
\]
and this is well-defined since it is a continuous linear operator
between Banach spaces.

\subsection{Kernel of trace operator on Lipschitz domain}

Now we are going to complete the proof by using change of variables
formulas into the proof of Theorem 2 in §5.5 of \cite{key-3}.

\paragraph{Theorem 4 (Theorem 4.10, §2.4.3 in \cite{key-2}). }

Let $D\subset\mathbb{R}^{d}$ be a Lipschitz domain and $u\in W^{1,p}(D)$,
$p\in[1,\infty)$. Then $u\in\overline{C_{c}^{\infty}\left(D\right)}$
if and only if $Tu=0$ on $\partial D$.

Proof: One side is trivial. To show the converse, firstly we establish
a priori estimate.\\
Let $Tu=0$ on $\partial D$ and $a\in\partial D$. Then there exist
$r>0$, Lipschitz map $\gamma$, $F:\,\mathbb{R}^{d}\rightarrow\mathbb{R}^{d}$,
$F(x)=y$ be Lipschitz coordinate transformation as in §1.1, and $u_{m}\in\mathcal{C}^{1}(D)$
such that $u_{m}\rightarrow u$ in $W^{1,p}(D)$ and $Tu_{m}\rightarrow0$
in $L^{p}(\partial D)$ as $m\rightarrow\infty$. If $y'\in\mathbb{R}^{d-1}$,
$y_{d}>0$, and $(y',y_{d})\in F(B(a,r)\cap\partial D)$, then
\begin{equation}
|u_{m}(y',y_{d})|\leq|u_{m}(y',0)|+\stackrel[0]{y_{d}}{\int}|\partial_{y_{d}}u_{m}(y',s)|\,ds.
\end{equation}
We use Theorem 2 (or just ordinary change variables) by taking $m=n=d-1,$
\[
g(y',y_{d})=|u_{m}(y',y_{d})|^{p}\mathtt{1}_{F(B(a,r)\cap\partial D)}(y'),
\]
 and $f(x')=y'$, thus $Jf=1$. Then
\[
\underset{B(a,r)\cap\partial D}{\int}|u_{m}(x',x{}_{d})|^{p}\,dx'=\underset{F(B(a,r)\cap\partial D)}{\int}|u_{m}(y',0)|^{p}\,d\mathcal{H}^{d-1}(y').
\]
Taking p power (so we need $p<\infty$) on equation (1), we have 
\[
|u_{m}(y',y_{d})|^{p}\leq C(|u_{m}(y',0)|^{p}+(\stackrel[0]{y_{d}}{\int}|\partial_{y_{d}}u_{m}(y',s)|\,ds)^{p}),
\]
and then
\[
(\stackrel[0]{y_{d}}{\int}|\partial_{y_{d}}u_{m}(y',s)|\,ds)^{p}\leq y_{d}^{p-1}\stackrel[0]{y_{d}}{\int}|\partial_{y_{d}}u_{m}(y',s)|^{p}\,ds
\]
by Jensen's inequality. Then integrate with respect to $y'$, on $B:=\{y'\in\mathbb{R}^{d-1}\mid(y',\cdot)\in F(B(a,r))\}$
\[
\underset{B}{\int}|u_{m}(y',y_{d})|^{p}\,dy'\leq C(\underset{B}{\int}|u_{m}(y',0)|^{p}\,dy'+y_{d}^{p-1}\stackrel[0]{y_{d}}{\int}\underset{B}{\int}|Du_{m}(y',s)|^{p}\,dy'ds).
\]
Let $m\rightarrow\infty$, and then we have a priori estimate
\begin{equation}
\underset{B}{\int}|u(y',y_{d})|^{p}\,dy'\leq Cy_{d}^{p-1}\stackrel[0]{y_{d}}{\int}\underset{B}{\int}|Du(y',s)|^{p}\,dy'ds,
\end{equation}
for a.e. $y_{d}>0$.\\
 Now we are going to approximate $u$ under Lipschitz coordinate transformation.\\
Let $\zeta\in\mathcal{C}_{c}^{\infty}(\mathbb{R}_{+})$ such that
$0\leq\zeta\leq1$, $\zeta|_{\mathbb{R}_{+}\setminus[0,2]}=0$, $\zeta|_{[0,1]}=1$,
$\zeta_{k}(y):=\zeta(ky_{d})$, $\forall y\in\mathbb{R}_{+}^{d}$,
and $w_{k}:=u(1-\zeta_{k})$$\rightarrow u$ in $L^{p}(\mathbb{R}_{+}^{d})$
as $k\rightarrow\infty$. Note that $\sup|\zeta'|<\infty$. The remainder
is to estimate
\[
|Dw_{k}-Du|^{p}=|\zeta_{k}Du+ku\zeta'|^{p}\leq C(|\zeta_{k}|^{p}|Du|^{p}+k^{p}|\zeta'(ky_{d})|^{p}|u|^{p}).
\]
Since $supp\zeta{}_{k}\subset[0,\frac{2}{k}]$, $\int|\zeta_{k}|^{p}|Du|^{p}\rightarrow0$
as $k\rightarrow\infty$ by Lebesgue's dominated convergence Theorem.
Since $supp\zeta'\subset[0,2]$, to integrate the last term on $F^{-1}(B\times[0,2])$,
we use Theorem 2 again by taking $m=n=d,$ 
\[
g(y',y_{d})=|u(y',y_{d})|^{p}\mathtt{1}_{B}(y')\mathtt{1}_{[0,2/k]}(y_{d}),
\]
 and $F(x',x_{d})=(y',y_{d})$, thus $JF=1$ a.e. And then the priori
estimate (2) shows
\begin{equation}
Ck^{p}\stackrel[0]{2/k}{\int}\underset{B}{\int}|u(y',y_{d})|^{p}\,dy'dy_{d}\leq Ck^{p}(\stackrel[0]{2/k}{\int}y_{d}^{p-1}\,dy_{d})(\stackrel[0]{2/k}{\int}\underset{B}{\int}|Du(y',s)|^{p}\,dy'ds),
\end{equation}
therefore
\[
(3)\leq\stackrel[0]{2/k}{\int}\underset{B}{\int}|Du(y',s)|^{p}\,dy'ds\rightarrow0
\]
as $k\rightarrow\infty$. By partition of unity as in §1.1 above,
we deduce $w_{k}\rightarrow u$ in $W^{1,p}(\mathbb{R}_{+}^{d})$,
and $w_{k}=0$ if $0<y_{d}<1/k$. To conclude, we can mollify $w_{k}$
to produce $u_{k}\in\mathcal{C}_{c}^{\infty}(\mathbb{R}_{+}^{d})$
such that $u_{k}\rightarrow u$ in $W^{1,p}(\mathbb{R}_{+}^{d})$,
i.e., $u\in\overline{C_{c}^{\infty}\left(\mathbb{R}_{+}^{d}\right)}$.
This completes the proof.

\section{Concluding Remark}

The only difference between the Theorem 2 in §5.5 of \cite{key-3},
and this Theorem 4, is that we use change of variables formulas for
Lipschitz coordinate transformation. Recall that we need a topological
condition: boundary of Lipschitz domain $\partial D$ is compact,
to make sure that the trace operator is bounded. However we have no
idea what would happen, if $D$ is not bounded, or $p=\infty$, even
in smooth domain.\\
\\

\lyxaddress{e-mail: 81040001s@ntnu.edu.tw}

\end{document}